\documentclass[12pt]{article}

\usepackage{amsmath}
\usepackage{amsfonts}
\usepackage{amssymb}
\usepackage{graphicx}
\usepackage[english]{babel}
\usepackage{amsthm}
\usepackage{geometry}
\usepackage{tabto}
\usepackage{latexsym}
\usepackage{enumitem}
\usepackage{float}
\usepackage{caption}
\usepackage{subcaption}
\usepackage{multicol}
\usepackage[utf8]{inputenc}
\usepackage[english]{babel}
\usepackage{setspace}
\usepackage[toc,page]{appendix}
\usepackage{mathtools}

\newtheorem{theorem}{Theorem}

\newtheorem{lemma}{Lemma}
\newtheorem{proposition}{Proposition}
\newtheorem{definition}{Definition}

\usepackage{authblk}
\usepackage{varwidth}

\providecommand{\keywords}[1]{\textbf{\textit{Keywords---}} #1}

\title{On Finding the Eigenvalues of the Matrix of Rotation Symmetric Boolean Functions}
\author{Manuel Albrizzio}
\affil{\protect\begin{varwidth}[t]{\linewidth}\protect\centering Department of Mathematics, Bucknell University, One Dent Drive, \par Lewisburg, PA  17837 \protect\end{varwidth}} 
\date{}

\begin{document}
	
	\maketitle
	
	\begin{abstract}
		 We consider the action on $\mathbb{F}_2^n$ by cyclic permutations ($\mathbb{Z}/n\mathbb{Z}$). Two elements $x, y\in \mathbb{F}_2^n$ are in the same orbit if they are cyclic shifts of each other. Cryptographic properties of rotation symmetric Boolean functions can be efficiently computed using the square matrix $_n\mathcal{A}$, the construction of which uses orbit representatives of the cyclic shifting action. In 2018, Ciungu and Iovanov proved that $_n\mathcal{A}^2=2^n\cdot I$, the identity matrix of dimension $g_n\times g_n$ where $g_n$ is the number of orbits. In this paper, we answer the open question of the precise number of positive and negative eigenvalues of $_n\mathcal{A}$. 
	\end{abstract}
	
	
	\begin{section}{Introduction}
		Introduced by Pieprzyk and Qu in their 1998 paper \cite{Pieprzyk}, Rotation Symmetric Boolean Functions (RSBFs) serve as important components in hashing algorithms. They allow for fast and efficient evaluations in schemes involving Boolean functions. A Boolean function is a function $f: \mathbb{F}_2^n \rightarrow \mathbb{F}_2$ where $\mathbb{F}_2=\{0, 1\}$, and such a function is an RSBF if $f(x_1,x_2,\dots, x_n) = f(x_2,x_3,\dots, x_n,x_1)$, i.e. the function remains invariant under cylic shifting of an element of $\mathbb{F}_2^n$. The space $\mathbb{F}_2^n$ can be partitioned by this cyclic shifting, an equivalence relation, whose equivalence classes we denote by $G_{n,i}$. Alternatively, $G_{n,i}$ are the orbits obtained by the natural action of the cyclic group $C_n$ on $\mathbb{F}_2^n$. Let $g_n$ be the number of such orbits and let $\{\Lambda_{n,i}\}$, $ \Lambda_{n,i}\in G_{n,i}$, be the orbit representatives which come first lexicographically in each orbit. A useful tool when computing the security and cryptographic properties of RSBFs, which is also considered by several authors \cite{Stanica,Maximov}, is the matrix
		
		$$ _n\mathcal{A} = \left( \displaystyle\sum_{x\in G_{n,i}} (-1)^{x \cdot \Lambda_{n,j}} \right)_{i,j} .$$
		
			Ciungu and Iovanov showed \cite{Ciungu} that the square of this matrix is a multiple of the $g_n \times g_n$ identity matrix; specifically
		
		$$ _n\mathcal{A} = 2^n \cdot \text{Id}_{g_n} .$$
		
			This paper answers the open question posed in the same paper; they asked to precisely determine the number of positive and negative eigenvalues of $_n\mathcal{A}$. While this question looks trivial given the relation above, a closer analysis is needed of the case when $n$ is even.
			
			In this paper, we look first at some preliminary definitions about Boolean functions and formally defining RSBFs. We then introduce the matrix and its construction. The last section contains proofs of lemmas and propositions leading up to the proof of our main theorem. 
			
	\end{section}
	
	\subsection*{Acknowledgements}
	This paper would not have been possible with the support and guidance of my advisors, Lavinia Ciungu and Miodrag Iovanov. I am extremely grateful for this. Additionally, I thank Nathan Ryan at Bucknell University for knowledge he has provided.

	\begin{section}{Preliminaries}
		
		Let $\mathbb{F}_2^n$ be the vector space of dimension $n$ over the two-element field $\mathbb{F}_2$. For two vectors $x=(x_1,x_2,\dots,x_n)$ and $y=(y_1,y_2,
		\dots,y_n)$ in $\mathbb{F}_2^n$, we say $x \cdot y=x_1y_1 +x_2y_2+\dots +x_ny_n\in \mathbb{F}_2$ is the scalar product where the multiplication and addition are over $\mathbb{F}_2$.
		
		\begin{definition}
			A \textbf{Boolean function} $f$ is map from $\mathbb{F}_2^n$ to $\mathbb{F}_2$.
		\end{definition}
	
		\begin{definition}
			For any $x\in \mathbb{F}_2^n$, the \textbf{Hamming weight} of $x$, denoted $wt(x)$, is the number of $1$'s in $x$. For any two elements $x$ and $y$ in $\mathbb{F}_2^n$, the \textbf{distance} between them, denoted $d(x,y)$, is the number of components the two elements differ from each other. We can see $d(x,y)=wt(x\oplus y)$ where $\oplus$ is the addition defined on $\mathbb{F}_2^n$.
		\end{definition}
	
		 Hash functions are essential in the exchange of long messages by securely creating a fixed-length ``fingerprint''  of the message that is used as a digital signature. Digital signatures provide a message with integrity, authentication, and nonrepudiation. The addition of a hash function gives the message an extra layer of security. Pieprzyk and Qu \cite{Pieprzyk} introduced the use of functions that are rotationally symmetric in iterations of hashing algorithms for efficient evaluations.
		 
   		Consider the linear function $\rho_n: \mathbb{F}_2^n\rightarrow \mathbb{F}_2^n$ that shifts components of the vector cyclically, i.e for $x=(x_1,x_2,\dots, x_n)\in \mathbb{F}_2^n$, $\rho_n(x)=(x_2,\dots,x_n,x_1)$. 
   		
   		\begin{definition}
   			A Boolean function $f$ is \textbf{rotation symmetric} if for any element $x \in \mathbb{F}_2^n$,
   			$$ f(\rho_n^k(x))=f(x)$$
   			for all $k$, $1\leq k \leq n$. 
   		\end{definition}
   	
   	We fix $n>2$. Let $C_n$ be the cyclic group generated by $\rho_n$, which we think of as a subgroup of the symmetric group on $n$ elements. Then a function $f$ is rotation symmetric if and only if it takes the same value on each orbit of the action of $C_n$ on $\mathbb{F}_2^n$. We denote $G_n(x)$ as the orbit of $x$ under this action, i.e $G_n(x)=G_n((x_1,x_2,\dots,x_n))=\{\rho^k(x)|1\leq k\leq n\}$. Let $g_n$ denote the number of orbits. A direct application of Burnside's Lemma shows
   	\begin{equation}
   		g_n=\dfrac{1}{n} \displaystyle\sum_{k|n} \phi(k)\cdot 2^{\frac{n}{k}}
   		\label{orbitnum} ,
   	\end{equation}
   	where $\phi(k)$ denotes the Euler's $\phi$-function. In \cite{Ciungu}, it is shown $g_n$ is even for all $n>2$.
   	
   	Lemma 1 in \cite{Stanica} shows the Walsh transform is constant on the orbits. Then the \textit{Walsh spectrum}, the vector of possible Walsh tranforms, is at most $g_n$ valued. \\
   	
   	\begin{definition}
   		The \textbf{Walsh transform}  of a Boolean function is a map $W_f:\mathbb{F}_2^n\rightarrow \mathbb{R}$, defined by
   		$$ W_f(w) = \displaystyle\sum_{x\in \mathbb{F}_2^n} (-1)^{f(x)+ x\cdot w} .$$
   	\end{definition} \vspace{1cm}
   	
   	\noindent Many cryptographic properties, such as non-linearity and resilience, of Boolean functions can be described using the Walsh spectrum.
   	
   	\begin{definition}
   		A Boolean function $f$ is considered \textbf{bent} if the Walsh transform has constant absolute value.
   	\end{definition}
		
	\end{section}

	\begin{section}{The matrix of rotation symmetric Boolean functions}
		
		Consider the orbit representatives $\Lambda_{n,i}$ $(0\leq i \leq g_n-1)$ which are the elements that come first in lexicographical order in each orbit. Using the set $\{\Lambda_{n,i}\}$ (ordered lexicographically as well), we construct the matrix related to the set of $n$-variables RSBFs by:
		
		$$ _n\mathcal{A} = \left( \displaystyle\sum_{x\in G_{n,i}} (-1)^{x \cdot \Lambda_{n,j}} \right)_{i,j} .$$
		
		Dalai \cite{Dalai} studied a variation of this matrix $_n\mathcal{A}^{\pi}$ where the columns were permuted. For our purposes, however, we stick with the convention in \cite{Ciungu}. \\
		
		\textbf{Example.}  Let $n=4$. The orbits of $\{0,1\}^n$ are as follows: \\
		$G_4(0,0,0,0)=\{(0,0,0,0)\},$ \\
		$G_4(0,0,0,1)=\{(0,0,0,1),(0,0,1,0),(0,1,0,0),(1,0,0,0)\},$ \\
		$G_4(0,0,1,1)=\{(0,0,1,1),(0,1,1,0),(1,1,0,0),(1,0,0,1)\},$ \\
		$G_4(0,1,0,1)=\{(0,1,0,1),(1,0,1,0)\},$ \\
		$G_4(0,1,1,1)=\{(0,1,1,1),(1,1,1,0),(1,1,0,1),(1,0,1,1)\},$ \\
		$G_4(1,1,1,1)=\{(1,1,1,1)\}.$
		
		We can see the number of partitions is $g_n=6$, matching \eqref{orbitnum}. The representatives should be chosen to be $\Lambda_{4,0}=(0,0,0,0)$, $\Lambda_{4,1}=(0,0,0,1)$, $\Lambda_{4,2}=(0,0,1,1)$, $\Lambda_{4,3}=(0,1,0,1)$, $\Lambda_{4,4}=(0,1,1,1)$, and $\Lambda_{4,5}=(1,1,1,1)$, as they are the first in lexicographic order. The resulting matrix is then
		
		$$_n\mathcal{A}=\begin{pmatrix}
			1 & 1 & 1 & 1 & 1 & 1 \\
			4 & 2 & 0 & 0 & -2 & -4 \\
			4 & 0 & 0 & -4 & 0 & 4 \\
			2 & 0 & -2 & 2 & 0 & 2 \\
			4 & -2 & 0 & 0 & 2 & -4 \\
			1 & -1 & 1 & 1 & -1 & 1
		\end{pmatrix} .$$ \\
	
	We now state the main result in \cite{Ciungu} about this matrix:
	
	\begin{theorem}[\cite{Ciungu}]
		$_n\mathcal{A}^2=2^n\cdot I$, where I denotes the identity matrix of size $g_n\times g_n$.
		\label{squareid}
	\end{theorem}
	
	The above theorem was proved in \cite{Ciungu} using a direct approach and also using characters of the semidirect product of the group $\mathbb{F}_2^n \rtimes (\mathbb{Z}/n\mathbb{Z})$. 
	
	We now state the main theorem of this paper.
	
	\begin{theorem}
		For $n>2$, the matrix $_n\mathcal{A}$ has eigenvalues $\pm 2^{n/2}$ that satisfy the following:
		\begin{itemize}
			\item When $n$ is odd, both eigenvalues appear $g_n/2$ times.
			
			\item When $n$ is even,
				$$\# \text{ positive eigenvalues }= \dfrac{g_n}{2} + \dfrac{1}{2n}\displaystyle\sum_{\substack{k|n \\ 2|k}} \phi\left(k\right) 2^{n/k} ,$$
				$$\# \text{ negative eigenvalues }= \dfrac{g_n}{2} - \dfrac{1}{2n}\displaystyle\sum_{\substack{k|n \\ 2|k}} \phi\left(k\right) 2^{n/k} .$$
		\end{itemize}
		\label{Mainthm}
	\end{theorem}
	
	The proof of this is very involved. We present it as a series of lemmas and propositions leading to our final result.
	
	\begin{lemma}
		For odd $n>2$, the matrix $_n\mathcal{A}$ has eigenvalues $\pm 2^{n/2}$ with multiplicity $\frac{g_n}{2}$ for each.
		\label{oddn}
	\end{lemma}
	
	\textbf{Proof:} From Theorem \ref{squareid}, we see the matrix $_n\mathcal{A}$ satisfies the equation $p(x)=x^2-2^n$, which we consider as an element of the polynomial ring over the rationals, $\mathbb{Q}[x]$. As $n$ is odd, $p(x)$ is irreducible over the rationals. As $\mathbb{Q}[x]$ is a PID,  this shows $p(x)$ is the minimal polynomial of the matrix and also its only invariant factor. The rational canonical form of $_n\mathcal{A}$ is then a block matrix with the following companion matrix in each block
	$$\begin{bmatrix}
		0 & 2^n \\
		1 & 0
	\end{bmatrix}.$$
	
	This matrix has eigenvalues $\pm 2^{n/2}$ each with multiplicity 1. Therefore, the eigenvalues $\pm 2^{n/2}$ appear with equal multiplicity. \qed \\
	
	Unlike in the proof of the odd case, the case when $n$ is even gives us a reducible minimal polynomial (which is still $p(x)$). Hence the invariant factors could contain copies of either $x-2^{n/2}$ or $x+2^{n/2}$. This actual factor and the number of copies appearing is unclear. We focus on the trace of the matrix to prove the even case. From the construction, we see the trace is given by.
	\begin{equation}
		\mathrm{Tr}(_n\mathcal{A})=\displaystyle\sum_{i=1}^{g_n} \displaystyle\sum_{x\in G_n(\Lambda_{n,i})} (-1)^{x \cdot \Lambda_{n,i}},
		\label{traceog}
	\end{equation}
	
	\noindent where $\Lambda_{n,i}$ are the representative elements of the orbits of $\mathbb{F}_{2}^n$ when acted on by $C_n=\langle\rho_n\rangle$, the group generated by a cyclic permutation.
	
	\begin{proposition}
		$$\mathrm{Tr}\left(_n\mathcal{A}\right) = \dfrac{1}{n} \displaystyle\sum_{\sigma\in C_n} \displaystyle\sum_{x\in \mathbb{F}_2^n} (-1)^{x\cdot \sigma x}$$
		\label{tracenew}
	\end{proposition}
	
	\textbf{Proof:} Each $x\in G_n(\Lambda_{n,i})$ in \eqref{traceog} is a shift of $\Lambda_{n,i}$ by some number of spaces for each $i$. Using the orbit-stabilizer relation gives us
	
	\begin{equation}
		\mathrm{Tr}(_n\mathcal{A})=\displaystyle\sum_{i=1}^{g_n} \left(\displaystyle\sum_{\sigma\in C_n} (-1)^{\Lambda_{n,i}\cdot \sigma \Lambda_{n,i}}\right) \cdot \frac{|G_n(\Lambda_{n,i})|}{n}
	\end{equation}
	
	Let $\tau\in C_n$ and $y,x\in \mathbb{F}_{2^n}$. Then, if $y=\tau x$, we know
	
	\begin{equation}
		\displaystyle\sum_{\sigma \in C_n} (-1)^{y\cdot \sigma y} = \displaystyle\sum_{\sigma \in C_n} (-1)^{x\cdot \sigma x},
		\label{equalsum}
	\end{equation}
	
	\noindent since $\tau x\cdot \sigma \tau x=\tau x\cdot \tau \sigma x=x\cdot \sigma x$. Writing the cardinality of the orbit as a sum, \eqref{equalsum}  yields 
	
	\begin{equation}
		\mathrm{Tr}(_n\mathcal{A})= \dfrac{1}{n}\displaystyle\sum_{i=1}^{g_n} \displaystyle\sum_{x\in (G_n(\Lambda_{n,i}))} \displaystyle\sum_{\sigma\in C_n} (-1)^{x\cdot \sigma x} = \dfrac{1}{n}\displaystyle\sum_{x\in \mathbb{F}_{2}^n} \displaystyle\sum_{\sigma\in C_n} (-1)^{x\cdot \sigma x} .
	\end{equation}
	
	Switching the order of summation proves our proposition. \qed \\
	
	We now claim the following value of the inner sum in Proposition \ref{tracenew}: For a given $\sigma\in C_n$,
	
	\begin{equation}
		S_\sigma:=\displaystyle\sum_{x\in \mathbb{F}_2^n} (-1)^{x\cdot \sigma x} =
		\begin{cases}
			2^{\frac{n}{2}+\frac{n}{ord(\sigma)}}    &   ord(\sigma) \text{ is even }\\
			0                                  &    ord(\sigma) \text{ is odd }
		\end{cases}
		\label{innsum}
	\end{equation}
	
	The proof of \eqref{innsum} varies slightly depending on the order of the shift $\sigma$, namely whether or not $ord(\sigma)$ divides $n$ or not. Let us consider the case when $ord(\sigma)$ divides $n$.
	
	\begin{lemma}
		Let $\sigma\in C_n$ and $k$ be a positive integer such that $k\cdot ord(\sigma)=n$. Then $S_\sigma$ takes the value as in \eqref{innsum}.
		\label{innsumdivide}
	\end{lemma}
	
	 \noindent\textbf{Proof:} Let $\sigma=\rho^k$ where $k|n$. For any $v\in \mathbb{F}_{2^n}$, we divide the vector in blocks of size $k$:
	$$v=(a_0, a_1, \dots, a_{n-1})= (\textbf{b}_1, \textbf{b}_2, \dots, \textbf{b}_{n/k}) \implies \sigma v=(\textbf{b}_2, \dots, \textbf{b}_{n/k}, \textbf{b}_1).$$
	Then
	$$v\cdot \sigma v=\textbf{b}_1\cdot \textbf{b}_2 + \textbf{b}_2\cdot \textbf{b}_3 +\dots +\textbf{b}_{n/k}\cdot \textbf{b}_1.$$
	We now rewrite the sum in the following way
	
	\begin{equation}
		\displaystyle\sum_{v\in \mathbb{F}_2^n} (-1)^{v\cdot \sigma v}=\displaystyle\sum_{\substack{\textbf{b}_i\in \mathbb{F}_{2}^k \\ 1\leq i\leq n}} (-1)^{\textbf{b}_1\cdot \textbf{b}_2 + \textbf{b}_2\cdot \textbf{b}_3 +\dots +\textbf{b}_{n/k}\cdot \textbf{b}_1}.
		\label{dotprodeqsum}
	\end{equation}
	
	By the first isomorphism theorem, the multiplicative map $\phi_a:\mathbb{F}_{2}^n\rightarrow \mathbb{F}_2$ mapping $x\mapsto a\cdot x$ for $a\in \mathbb{F}_{2}^n$ has a kernel of size $2^n$ when $a=\mathbf{0}$ or size $2^{n-1}$ otherwise. This means the sum
	$ \displaystyle\sum_{c\in \mathbb{F}_2^n} (-1)^{a\cdot c} $
	is either $2^n$ if $a=\textbf{0}$, or 0 otherwise since $a\cdot c$ equals 0 and 1 the same number of times.
	
	Let us separate our analysis into two cases:
	
	\begin{itemize}
		
		\item Assume $\frac{n}{k}$ is odd. The exponent in \eqref{dotprodeqsum} then has an odd number of terms. This gives us
		
		\begin{equation}
			\resizebox{.865\hsize}{!}{$\begin{aligned}
				\displaystyle\sum_{v\in \mathbb{F}_2^n} (-1)^{v\cdot \sigma v} &=
				\displaystyle\sum_{\substack{\textbf{b}_i\in \mathbb{F}_{2}^k \\ 1\leq i\leq n}} (-1)^{\textbf{b}_1\cdot \textbf{b}_2 + \textbf{b}_2\cdot \textbf{b}_3 +\dots +\textbf{b}_{n/k}\cdot \textbf{b}_1} \\
				&= \displaystyle\sum_{\substack{\textbf{b}_i\in \mathbb{F}_{2}^k \\ 1\leq i\leq n \\ i \text{ is odd }}} (-1)^{\textbf{b}_1 \textbf{b}_{n/k}}\left[  \displaystyle\sum_{\substack{\textbf{b}_j\in \mathbb{F}_{2}^k \\ 1\leq j\leq n \\ j \text{ is even }}} (-1)^{\textbf{b}_2\cdot(\textbf{b}_1\oplus\textbf{b}_3)+\textbf{b}_4\cdot(\textbf{b}_3\oplus\textbf{b}_5)+\dots+\textbf{b}_{n/k-1}\cdot(\textbf{b}_{n/k-2}\oplus\textbf{b}_{n/k})}\right] \\
				&= \displaystyle\sum_{\substack{\textbf{b}_i\in \mathbb{F}_{2}^k \\ 1\leq i\leq n \\ i \text{ is odd }}} (-1)^{\textbf{b}_1 \textbf{b}_{n/k}} \left[\displaystyle\sum_{\substack{\textbf{b}_j\in \mathbb{F}_{2}^k \\ 1\leq j\leq n \\ j \text{ is even }}} (-1)^{\textbf{b}_2\cdot(\textbf{b}_1\oplus\textbf{b}_3)}(-1)^{\textbf{b}_4\cdot(\textbf{b}_3\oplus\textbf{b}_5)}\dots(-1)^{\textbf{b}_{n/k-1}\cdot(\textbf{b}_{n/k-2}\oplus\textbf{b}_{n/k})} \right]
			\end{aligned}$}
		\end{equation}

		Let us consider the first product in the inner sum of the above, $\sum\limits_{\textbf{b}_2\in \mathbb{F}_2^k} (-1)^{\textbf{b}_2\cdot (\textbf{b}_1\oplus\textbf{b}_3)}$. The term in the sum is exactly the map of dot-product multiplication by $\textbf{b}_2$. There are two cases to this product. One case is the term which will be 1 for $2^{k-1}$ vectors in $\mathbb{F}_2^k$ and $-1$ for the other $2^{k-1}$ vectors. This occurs when the term in the parentheses is not \textbf{0}. The value of the sum would then be 0. The other case is the term in the sum will be 1 for all vectors in $\mathbb{F}_2^k$. This occurs when the term in the parentheses is \textbf{0}. We know the only time when $\textbf{b}_1\oplus\textbf{b}_3=0$ is when $\textbf{b}_1=\textbf{b}_3$. We continue in this way; i.e for any fixed even $j$, $1\leq j\leq n$,
		$$\displaystyle\sum_{\textbf{b}_j\in \mathbb{F}_{2}^k} (-1)^{\textbf{b}_j\cdot (\textbf{b}_{j-1}\oplus\textbf{b}_{j+1})}=
		\begin{cases}
			2^k & \textbf{b}_{j-1}=\textbf{b}_{j+1}\\
			0 & \textbf{b}_{j-1} \neq \textbf{b}_{j+1}
		\end{cases} .
		$$
		There are $\frac{n/k -1}{2}$ products. This gives us
		$$\begin{aligned} 
			S_\sigma &= \displaystyle\sum_{\textbf{b}_1\in\mathbb{F}_{2}^k} (-1)^{\textbf{b}_1\cdot \textbf{b}_1}\cdot 2^{\frac{(n/k-1)\cdot k}{2}} \ \ \ \ \text{since} \ \textbf{b}_1=\textbf{b}_3=\dots=\textbf{b}_{n/k} \\ \\
			&=2^{\frac{(n/k-1)\cdot k}{2}} \displaystyle\sum_{\textbf{b}_1\in\mathbb{F}_{2}^k} (-1)^{wt(\textbf{b}_1)} \\ \\
			&=2^{\frac{(n/k-1)\cdot k}{2}}\left[ \binom{k}{0} - \binom{k}{1} + \dots (-1)^k \binom{k}{k} \right] \\ \\
			&=0.\end{aligned}$$
		The above uses the fact that
		$$\displaystyle\sum_{i=0}^n (-1)^i\binom{n}{i} =(-1+1)^n=0.$$
		
		\item Now assume $\frac{n}{k}$ is an even number. The exponent in \eqref{dotprodeqsum} has an even number of terms. Hence,
		$$\displaystyle\sum_{\substack{\textbf{b}_i\in \mathbb{F}_{2}^k \\ 1\leq i\leq n}} (-1)^{\textbf{b}_1\cdot \textbf{b}_2 + \textbf{b}_2\cdot \textbf{b}_3 +\dots +\textbf{b}_{n/k}\cdot \textbf{b}_1} $$
		$$= \displaystyle\sum_{\substack{\textbf{b}_i\in \mathbb{F}_{2}^k \\ 1\leq i\leq n \\ i \text{ is odd }}} \displaystyle\sum_{\substack{\textbf{b}_j\in \mathbb{F}_{2}^k \\ 1\leq j\leq n \\ j \text{ is even }}} (-1)^{\textbf{b}_2\cdot(\textbf{b}_1\oplus\textbf{b}_3)+\textbf{b}_4\cdot(\textbf{b}_3\oplus\textbf{b}_5)+\dots+\textbf{b}_{n/k}\cdot(\textbf{b}_{n/k-1}\oplus\textbf{b}_{1})}.$$
		The term in the sum can be written as follows:
		
		\begin{equation}
			\displaystyle\sum_{\substack{\textbf{b}_i\in \mathbb{F}_{2}^k \\ 1\leq i\leq n \\ i \text{ is odd }}} \displaystyle\sum_{\substack{\textbf{b}_j\in \mathbb{F}_{2}^k \\ 1\leq j\leq n \\ j \text{ is even }}} (-1)^{\textbf{b}_2\cdot(\textbf{b}_1\oplus\textbf{b}_3)}(-1)^{\textbf{b}_4\cdot(\textbf{b}_3\oplus\textbf{b}_5)}\dots(-1)^{\textbf{b}_{n/k}\cdot(\textbf{b}_{n/k-1}\oplus\textbf{b}_{1})}.
			\label{sumproduct}
		\end{equation}
		
		As in the previous case, we run into a sum of terms of the dot-product multiplication map. Let us apply the result of the value here, i.e for any fixed even $j, \  1\leq j\leq n$
		$$\displaystyle\sum_{\textbf{b}_j\in \mathbb{F}_{2}^k} (-1)^{\textbf{b}_j\cdot (\textbf{b}_{j-1}\oplus\textbf{b}_{j+1})}=
		\begin{cases}
			2^k & \textbf{b}_{j-1}=\textbf{b}_{j+1}\\
			0 & \textbf{b}_{j-1} \neq \textbf{b}_{j+1}
		\end{cases}.
		$$
		
		Since there is no product outside the inner sum as seen in the previous case, our value $S_\sigma$ is now nonzero. In fact,
		
		$$S_\sigma=2^{\frac{n/k\cdot k}{2}}\displaystyle\sum_{\textbf{b}_1\in \mathbb{F}_{2}^k} 1= 2^{\frac{n}{2}+k} $$
		as there were $\frac{n/k}{2}$ terms in \eqref{sumproduct} each equaling $2^k$.

		Given $\mathrm{ord}(\rho^k)=\frac{n}{k}$, our claim is proven. \qed

	\end{itemize}

	\begin{proposition}
		Let $\sigma\in C_n$. Then $S_\sigma$ takes the value as in \eqref{innsum}.
		\label{claiminnsum}
	\end{proposition}
	 
	 \noindent\textbf{Proof:} Let $\sigma=\rho^k$ where $ord(\rho^k)=d$, $k\nmid n$, and $\frac{n}{d}=\gcd(n,k)=m$. Given $v=(a_0,a_1,\dots, a_{n-1})\in \mathbb{F}_{2}^n$, we set the following vectors.
	
	$$\begin{aligned}
		\textbf{b}_1&=(a_0, a_1, \dots, a_{k-1}) \\
		\textbf{b}_2&=(a_k, a_{k+1}, \dots, a_{2k-1}) \\
		\vdots \\
		\textbf{b}_i&=(a_{(i-1)k}, a_{(i-1)k+1},\dots, a_{ik-1}) \\
		\vdots \\
		\textbf{b}_d&=(a_{(d-1)k}, a_{(d-1)k+1},\dots, a_{dk-1}) 
	\end{aligned} $$
	where the indices are$\mod n$. Its easy to see
	$$v\cdot \sigma v=m^{-1}\left[ \textbf{b}_1\cdot \textbf{b}_2 + \textbf{b}_2\cdot \textbf{b}_3 +\dots +\textbf{b}_{d-1}\cdot \textbf{b}_d  +\textbf{b}_d\cdot \textbf{b}_1\right] \in \mathbb{F}_2.$$
	
	There are $dk$ elements contained in the $d$ vectors, which means some components of $v$ will repeat in indices. Hence, we first will try to find the components that uniquely determine $v$, which will determine the other components in the $\textbf{b}_i$'s. Consider the following vectors
	$$\begin{aligned}
		\textbf{c}_1&=(a_0, a_1, \dots, a_{m-1}) \\
		\textbf{c}_2&=(a_k, a_{k+1}, \dots, a_{2m-1}) \\
		\vdots \\
		\textbf{c}_i&=(a_{(i-1)k}, a_{(i-1)k+1},\dots, a_{im-1}) \\
		\vdots \\
		\textbf{c}_d&=(a_{(d-1)k}, a_{(d-1)k+1},\dots, a_{dm-1}) ,
	\end{aligned} $$
	which are the first $m$ components of the $\textbf{b}$ vectors (in order). There are $dm=n$ components here, and in fact the indices ($\mod n$) are unique. Indeed if any two components were equivalent, there would exist a $1\leq i,j\leq d$ and $0\leq s,t \leq m-1$ such that 
	$$(i-1)k+s\equiv (j-1)k+t \mod n \implies (i-j)k+(s-t)\equiv 0\mod n.$$
	Then there exists an $h$ such that $(i-j)k+(s-t)=hn$. This means $m|(s-t)$, and that can only happen if $s=t$ since $|s-t|<m$. Hence, $i=j$ and two vectors are repeated. This is impossible as the order of $\sigma$ is $d$. Finally, we can see, based on where the vectors are created, that
	$$v\cdot \sigma v= \textbf{c}_1\cdot \textbf{c}_2 + \textbf{c}_2\cdot \textbf{c}_3 +\dots +\textbf{c}_{d-1}\cdot \textbf{c}_d  +\textbf{c}_d\cdot \textbf{c}_1 $$
	where each vector is in $\mathbb{F}_{2}^m$. By the same argument in Lemma \ref{innsumdivide}, we can conclude.
	$$S=
	\begin{cases}
		2^{\frac{d}{2}\cdot m +m}    &   d \text{ is even }\\
		0                                  &    d \text{ is odd }
	\end{cases}
	=
	\begin{cases}
		2^{\frac{n}{2} + \frac{n}{ord(\sigma)}}    &   ord(\sigma) \text{ is even }\\
		0                                  &    ord(\sigma) \text{ is odd }
	\end{cases}.
	$$
	This completes the proof of our claim. \qed
	
	We now have the parts to prove our theorem.
	
	\noindent\textbf{Proof of Theorem \ref{Mainthm}} When $n$ is odd, Lemma \ref{oddn} proves the result. Now let $n$ be even. From Proposition \ref{claiminnsum}, we know
	$$\mathrm{Tr}(_n\mathcal{A})=\dfrac{1}{n}\displaystyle\sum_{\substack{\sigma\in C_n \\ 2|ord(\sigma)}} 2^{\frac{n}{2} + \frac{n}{ord(\sigma)}} = \dfrac{1}{n} 2^{\frac{n}{2}} \displaystyle\sum_{\substack{1\leq k \leq n \\ 2|ord(\rho^k)}} 2^{\gcd(n,k)}.$$
	
	Assume $\gcd(n,k)=d$. Then $\#\{k : 1\leq k \leq n, \gcd(k,n)=d\} = \#\{ k: 1 \leq k \leq n, \gcd(k,\frac{n}{d})=1 \}=\phi(\frac{n}{d})$ where $\phi$ is the Euler's $\phi$-function. Using this, the trace formula looks similar to the formula for $g_n$:
	$$\mathrm{Tr}(_n\mathcal{A})=\dfrac{1}{n}2^{\frac{n}{2}}\displaystyle\sum_{\substack{d|n \\ 2|\frac{n}{d}}} \phi\left(\frac{n}{d}\right) 2^{d} .$$
	
	The trace is then always positive for even values of $n$ as all terms in the sum are positive. This means $_n\mathcal{A}$ has eigenvalue $2^{n/2}$ with multiplicity greater than or equal to the multiplicity of eigenvalue $-2^{n/2}$. Specifically, the matrix has $\frac{1}{n}\displaystyle\sum_{\substack{k|n \\ 2|k}} \phi\left(k\right) 2^{\frac{n}{k}} $ more positive eigenvalues, after rewriting. From this we can see
	
	$$\begin{aligned}\# \text{ positive eigenvalues } &= \dfrac{g_n}{2} + \frac{1}{2n}\displaystyle\sum_{\substack{k|n \\ 2|k}} \phi\left(k\right) 2^{\frac{n}{k}}  \ \ \text{and} \\
	\# \text{ negative eigenvalues } &= \dfrac{g_n}{2} - \frac{1}{2n}\displaystyle\sum_{\substack{k|n \\ 2|k}} \phi\left(k\right) 2^{\frac{n}{k}} \end{aligned}.$$ \qed

	\end{section}
	
	Knowing the eigenvalues provide more information about the matrix. We may use this in fast computations of cryptographic properties of RSBFs. In \cite{Stanica}, we see for an RSBF $f$ the value of the Walsh transforms is 
	
	$$ W_f(\Lambda_{n,j})=\displaystyle\sum_{i=0}^{g_n-1} (-1)^{f(\Lambda_{n,i})}\mathcal{A}_{i,j} $$
	
	\noindent where $\Lambda_{n,j}$ is an orbit representation and $A_{i,j}$ is the $(i,j)$-th element of $_n\mathcal{A}$. Moreover, $f$ is bent if and only if $ W_f(\Lambda_{n,j})= \displaystyle\sum_{i=0}^{g_n-1} (-1)^{f(\Lambda_{n,i})}\mathcal{A}_{i,j} = \pm 2^{n/2}$ for each $0\leq j \leq g_n-1$.

	\pagebreak

\end{document}